\documentclass[preprint]{imsart}

\usepackage[english]{babel}
\usepackage[utf8]{inputenc}
\usepackage{calc}

\usepackage{amsmath}
\usepackage{amsfonts}
\usepackage{amssymb}
\usepackage{amsthm}
\usepackage[numbers]{natbib}
\usepackage{xspace}
\usepackage{bm}
\usepackage{enumerate}
\usepackage[colorlinks,citecolor=blue,urlcolor=blue]{hyperref}

\startlocaldefs

\newtheorem{definition}{Definition}[section]

\newtheorem{lemma}[definition]{Lemma}
\newtheorem{theorem}[definition]{Theorem}
\newtheorem{proposition}[definition]{Proposition}
\newtheorem{corollary}[definition]{Corollary}

\theoremstyle{definition}

\newtheorem{remark}[definition]{Remark}

\numberwithin{equation}{section}

\DeclareMathOperator{\spanv}{span}

\newcommand{\Ai}{\ensuremath{\mathcal A}\xspace}
\newcommand{\Bi}{\ensuremath{\mathcal B}\xspace}
\newcommand{\Ci}{\ensuremath{\mathcal C}\xspace}

\newcommand{\Fi}{\ensuremath{\mathcal F}\xspace}
\newcommand{\Gi}{\ensuremath{\mathcal G}\xspace}

\newcommand{\Ni}{\ensuremath{\mathcal N}\xspace}
\newcommand{\Pii}{\ensuremath{\mathcal P}\xspace}

\newcommand{\Ti}{\ensuremath{\mathcal T}\xspace}
\newcommand{\Vi}{\ensuremath{\mathcal V}\xspace}

\newcommand{\Aiu}{\ensuremath{\mathcal A(u)}\xspace}

\newcommand{\zset}{\ensuremath{{\emptyset'}}\xspace}
\newcommand{\vset}{\ensuremath{{\emptyset}}\xspace}

\newcommand{\Cp}{\ensuremath{{ \bm{\mathcal{A}_{C}}} }\xspace}
\newcommand{\Gc}[1]{\ensuremath{{ \mathcal{G}_{\hspace{-0.4pt}#1}^{\hspace{0.5pt}*} }}\xspace}

\newcommand{\pc}[1]{\ensuremath{{ p_{\hspace{-0.2pt}#1} }}\xspace}
\newcommand{\Pc}[1]{\ensuremath{{ P_{\hspace{-0.2pt}#1} }}\xspace}

\newcommand{\Cpi}[1]{\ensuremath{{ U_{\hspace{-1pt}C}^{\hspace{0.2pt}#1}} }\xspace}
\newcommand{\eCpi}[1]{ \varepsilon_{#1} }

\newcommand{\vXCp}[1]{\ensuremath{ \bm{\mathbf{X}_{#1}} }\xspace}
\newcommand{\vxCp}[1]{\ensuremath{ \bm{\mathbf{x}_{#1}} }\xspace}

\newcommand{\pth}[1]{(#1)}
\newcommand{\pthb}[1]{\bigl(#1\bigr)}
\newcommand{\pthB}[1]{\Bigl(#1\Bigr)}
\newcommand{\pthbb}[1]{\biggl(#1\biggr)}

\newcommand{\bkt}[1]{[#1]}
\newcommand{\bktb}[1]{\bigl[#1\bigr]}
\newcommand{\bktB}[1]{\Bigl[#1\Bigr]}
\newcommand{\bktbb}[1]{\biggl[#1\biggr]}
\newcommand{\bktBB}[1]{\Biggl[#1\Biggr]}

\newcommand{\brc}[1]{\{#1\}}

\newcommand{\brcbb}[1]{\biggl\{#1\biggr\}}

\newcommand{\dt}{\ensuremath{\textrm d}\xspace} 
\newcommand{\disy}{\Delta} 
\newcommand{\eqdef}{\overset{\mathrm{def}}{=}}
\newcommand{\eqd}{\overset{\mathrm{d}}{=}}

\newcommand{\scpr}[2]{\left\langle #1,#2 \right\rangle}

\newcommand{\ivof}[1]{\ensuremath{\left(#1\right]}}

\newcommand{\ivff}[1]{\ensuremath{\left[#1\right]}}

\newcommand{\abs}[1]{\lvert#1\rvert}

\newcommand{\norm}[1]{\lVert#1\rVert}
\newcommand{\normb}[1]{\bigl\lVert#1\bigr\rVert}

\newcommand{\normbb}[1]{\biggl\lVert#1\biggr\rVert}

\renewcommand{\Pr}{\ensuremath{\mathbb P}\xspace}

\newcommand{\esp}[2][]{\mathbb{E}#1\bkt{#2}}
\newcommand{\espb}[2][]{\mathbb{E}#1\bktb{\hspace{1pt}#2\hspace{1pt}}}
\newcommand{\espB}[2][]{\mathbb{E}#1\bktB{#2}}
\newcommand{\espbb}[2][]{\mathbb{E}#1\bktbb{#2}}

\newcommand{\espc}[3][]{\mathbb{E}#1\bkt{\hspace{1pt}#2\hspace{1.5pt}|\hspace{1.5pt}#3\hspace{1pt}}}

\newcommand{\espcB}[3][]{\mathbb{E}#1\bktB{#2\Bigm|#3}}

\newcommand{\R}{\ensuremath{\mathbf{R}}\xspace}

\newcommand{\N}{\ensuremath{\mathbf{N}}\xspace}

\newcommand{\indi}{\ensuremath{\mathbf{1}}\xspace}

\DeclareMathOperator{\cov}{Cov}

\endlocaldefs

\begin{document}

\begin{frontmatter}

\title{A set-indexed Ornstein-Uhlenbeck process}
\runtitle{A set-indexed Ornstein-Uhlenbeck process}

\author{\fnms{Paul} \snm{Balan\c{c}a}%
\thanksref{t1}%
\ead[label=e1]{paul.balanca@ecp.fr}%
}%
\and%
\author{\fnms{ Erick} \snm{Herbin}%
\ead[label=e2]{erick.herbin@ecp.fr}%
\ead[label=u1,url]{www.mas.ecp.fr/recherche/equipes/modelisation\_probabiliste}%
}

\address{\printead{e1,e2}\\ \printead{u1}\\[1em] 
  \'Ecole Centrale Paris, 
  Laboratoire MAS\\ 
  Grande Voie des Vignes
-  92295 Ch\^atenay-Malabry, France
}

\thankstext{t1}{This article is based on a chapter from the author's Ph.D. thesis, prepared under the supervision of the second author.}

\affiliation{\'Ecole Centrale Paris}
\runauthor{P. Balan\c{c}a and E. Herbin}

\begin{abstract}
The purpose of this article is a set-indexed extension of the well-known Ornstein-Uhlenbeck process. 
The first part is devoted to a stationary definition of the random field and ends up with the proof of a complete characterization by its $L^2$-continuity, stationarity and set-indexed Markov properties. 
This specific Markov transition system allows to define a general \emph{set-indexed Ornstein-Uhlenbeck (SIOU) process} with any initial probability measure. 
Finally, in the multiparameter case, the SIOU process is proved to admit a natural integral representation.
\end{abstract}



\end{frontmatter}

\section{Introduction}

The study of multiparameter processes goes back to the 70' and the theory developed for years covers multiple properties of random fields (we refer to the recent books \cite{Khoshnevisan(2002)} and \cite{Adler.Taylor(2007)} for a modern review). For instance, Cairoli and Walsh \cite{Cairoli.Walsh(1975),Walsh(1978),Walsh(1986)} have deeply investigated the extension of the martingale and stochastic integral theories to the two-parameter framework. 
A vast literature also concerns the Markovian aspects of random fields. Similarly to the case of martingales, different interesting Markov properties can be formalized for multiparameter processes. Among these, the most commonly studied ones are \emph{sharp-Markov} \cite{Levy(1945),Dalang.Walsh(1992)}, \emph{germ-Markov} \cite{McKean(1963),Pitt(1971),Kunsch(1979)} and \emph{$\ast$-Markov} \cite{Cairoli(1971),Korezlioglu.Lefort.ea(1981)} properties. 
We refer to \cite{Balanca(2012)} for a more complete description of these concepts. 
The study multiparameter processes is still a very active area of research, particularly the analysis of sample paths and geometric properties (see e.g. \cite{Ayache.Shieh.ea(2011), Dalang.Nualart.ea(2011), Khoshnevisan.Xiao(2005), Tudor.Xiao(2007), Xiao(2009)a}).

Set-indexed processes constitute a natural generalization of multiparameter stochastic processes and their local regularity have been considered in the Gaussian case since the early work of Dudley \cite{Dudley(1973)} (see also \cite{Adler.Feigin(1984), Alexander(1986), Bass.Pyke(1985)}). Extending the literature on random fields, several different subjects have been recently investigated, including set-indexed martingales \cite{Ivanoff.Merzbach(2000)}, set-indexed Markov \cite{Ivanoff.Merzbach(2000)a, Balan.Ivanoff(2002), Balanca(2012)} and L\'evy processes \cite{Adler.Feigin(1984), Bass.Pyke(1984), Herbin.Merzbach(2011)},  and set-indexed fractional Brownian motion \cite{Herbin.Merzbach(2006), Herbin.Merzbach(2009)}. Although the set-indexed formalism appears to be more abstract, it usually offers a simpler and more condensed way to express technical concepts of multiparameter processes. For instance, the present work intensively uses the \Ci-Markov property introduced and developed in \cite{Balanca(2012)}. In the latter, the Chapman-Kolmogorov equation related to transition probabilities turns out to be more easily expressed using the set-indexed formalism than the two-parameter framework.

In this paper, we follow the framework established by Ivanoff and Merzbach in the context of set-indexed martingales \cite{Ivanoff.Merzbach(2000)}.
An {\em indexing collection} $\mathcal{A}$ is constituted of compact subsets of a locally compact metric space $\mathcal{T}$ equipped with a Radon measure on the $\sigma$-field generated by $\mathcal{A}$.
\Aiu and \Ci respectively denote the class of finite unions of sets belonging to \Ai and the collection of increments $C = A\setminus B$, where $A\in\Ai$ and $B\in\Aiu$. Finally, \zset denotes the set $\cap_{A\in\Ai} A$, which usually plays a role equivalent to $0$ in $\R_+^N$.
In the present article, we suppose that the collection $\Ai$ and the measure $m$ satisfy the following assumptions:
\begin{enumerate}[\it (i)]
  \item \zset is a nonempty set and \Ai is closed under arbitrary intersections;
  \item \emph{Shape hypothesis}: for any $A, A_1,\dotsc,A_k\in\Ai$ with $A\subseteq\cup_{i=1}^k A_i$, there exists $i\in\brc{1,\dotsc,k}$ such that $A\subseteq A_i$;
	\item $m(\zset) = 0$ and $m$ is monotonically continuous on \Ai, i.e. for any increasing sequence $(A_n)_{n\in\N}$ in \Ai, 
	\[
    \lim_{n\rightarrow \infty} m(A_n) = m\pthb{ \overline{\cup_{k\in\N}A_k} }.
  \]
\end{enumerate}
For sake of readability, we restrict properties of \Ai to the strictly required ones in the sequel. The particular case of $\mathcal{A}=\{ [0,t];\;t\in\R_+^N \}$ shows that the set-indexed formalism extends the multiparameter setting. Another simple example satisfying \emph{Shape} can be constructed on the $\R^3$-unit sphere: $\Ai = \brc{A_{\theta,\varphi}; \theta\in [0,\pi) \text{ and } \varphi\in [0,2\pi)}$ where $A_{\theta,\varphi} = \brc{(1,\widehat\theta,\widehat\varphi) : \widehat\theta\in\ivff{0,\theta} \text{ and } \widehat\varphi\in\ivff{0,\varphi}}$.    We refer to \cite{Ivanoff.Merzbach(2000)} for a more complete definition of an indexing collection used in the general theory of set-indexed martingales.

We investigate the existence and properties of a set-indexed extension of the Ornstein-Uhlenbeck (OU) process, originaly introduced in \cite{Uhlenbeck.Ornstein(1930)} and then widely used in the literature to represent phenomena in physics, biology and finance (e.g. see \cite{Frank.Daffertshofer.Beek(2000), {Lansky.Rospars(1995)}, Masuda(2004)}).
A well-known integral representation of the real-parameter OU process $X=\{X_t;\;t\in\R_+\}$ is given by
\begin{equation}  \label{eq:ou_integral}
\forall t\in\R_+;\quad  X_t = X_0\ e^{-\lambda t} + \int_0^t \sigma\ e^{\lambda(s-t)} \,\dt W_s,
\end{equation}
where $\lambda$ and $\sigma$ are positive parameters and the initial distribution $\nu=\mathcal{L}(X_0)$ is independent of the Brownian motion $W$. Furthermore, $X$ is a Markov process characterized by the following transition densities, for all $t\in\R_+$ and $x,y\in\R$;
\begin{equation} \label{eq:ou_kernel}
  p_t(x;y) = \frac{1}{\sigma_t\sqrt{2\pi}} \exp\bktbb{-\frac{1}{2\sigma^2_t}\pthb{ y - x e^{-\lambda t} }^2 } \qquad\text{where } \sigma^2_t = \frac{\sigma^2}{2\lambda}\pthb{1-e^{-2\lambda t}}.
\end{equation}
Two particular cases of initial distribution will be of specific interest in the sequel:
\begin{enumerate}[ \bf 1.]
	\item If $\nu = \delta_x$, $x\in\R$, $X$ is a Gaussian process with the following mean and covariance, for all $s,t\in\R_+$,
	\begin{equation}  \label{eq:ou_deltax}
	  \esp[_x]{X_t} = x e^{-\lambda t}\qquad\text{and}\qquad\cov_x(X_s,X_t) = \frac{\sigma^2}{2\lambda}\pthb{ e^{-\lambda\abs{t-s}} - e^{-\lambda(t+s)} }.
	\end{equation}
	\item If $\nu \sim \Ni\pth{ 0,\frac{\sigma^2}{2\lambda} }$, $X$ is a \emph{stationary Ornstein-Uhlenbeck process}, i.e. a zero-mean Gaussian process such that
	\begin{equation}  \label{eq:ou_stationary}
	  \forall s,t\in\R_+;\quad \esp[_\nu]{X_s X_t} = \frac{\sigma^2}{2\lambda}e^{-\lambda\abs{t-s}}.
	\end{equation}
\end{enumerate}

Since a set-indexed extension of the OU process cannot be directly derived from the integral representation \eqref{eq:ou_integral}, we first focus on the stationary process described in \eqref{eq:ou_stationary}. A natural way to extend this covariance to the set-indexed framework is to substitute the absolute value $\abs{t-s}$ with $d(U,V)$ where $d$ is a distance defined on the elements of $\mathcal{A}$. Similarly to the case of the set-indexed fractional Brownian motion described in \cite{Herbin.Merzbach(2006)}, we consider the choice $d(U,V)=m(U\disy V)$, where $\disy$ denotes the symmetric difference.

In Section \ref{sec:siou}, we first define a \emph{stationary set-indexed Ornstein-Uhlenbeck} process (ssiOU) as a zero-mean Gaussian process $X = \brc{X_U; U\in\Ai}$ such that
\begin{equation} \label{eq:ssiou_cov}
  \forall U,V\in\Ai;\quad \esp{X_U X_V} = \frac{\sigma^2}{2\lambda} e^{-\lambda m(U\disy V)},
\end{equation}
where $\lambda$ and $\sigma$ are positive parameters. Stationarity and Markov properties of this set-indexed process are studied, and lead to the complete characterization proved in Theorem~\ref{th:ssiou_charac}. Then, using the Markov kernel obtained, we are able to introduce in Definition~\ref{def:siou} a general \emph{set-indexed Ornstein-Uhlenbeck process}, whose law is consistent with the covariance structure \eqref{eq:ou_deltax} in particular case of initial Dirac distributions.

Finally, in Section \ref{sec:mpou}, we prove that in the multiparameter case, the set-indexed Ornstein-Uhlenbeck has a natural integral representation which generalizes expression \eqref{eq:ou_integral}.

\section{A stationary set-indexed Ornstein-Uhlenbeck process} \label{sec:siou}

In this section, we define a set-indexed extension of the stationary Ornstein-Uhlenbeck process, defined by the Gaussian covariance structure (\ref{eq:ou_stationary}).

\subsection{Definition and first properties}
 
As a preliminary to the definition, we need to prove that the expected covariance function of the process is positive definite in the same way as Lemma 2.9 of \cite{Herbin.Merzbach(2006)}.
 
\begin{lemma}\label{lem:cov}
If $\mathcal{A}$ is an indexing collection, $m$ a Radon measure on the $\sigma$-field generated by $\mathcal{A}$ and $\lambda, \sigma$ positive constants, the function $\Gamma : \mathcal{A}\times\mathcal{A}\rightarrow\R$ defined by
\begin{align*}
\forall U,V\in\Ai;\quad \Gamma(U,V) = \frac{\sigma^2}{2\lambda} e^{-\lambda m(U\disy V)},
\end{align*}
is positive definite.
\end{lemma}

\begin{proof}
 Let $f_1,f_2,\ldots,f_k$ be in $L^2(m)$ and $u_1,u_2,\ldots,u_k$ be in $\R$. Let \Vi be the vector space $\Vi=\spanv\pth{f_1,\ldots,f_k}$. Since $f\mapsto e^{-\frac{1}{2}\norm{f}_{L^2(m)}^2}$ is positive definite, there exists a Gaussian vector $X$ on the finite-dimensional space $\Vi$ such that 
  \[
    \forall \lambda>0,\ \forall f\in\Vi;\quad \espb{e^{i\sqrt{2\lambda}\scpr{X}{f}}} = e^{-\lambda\norm{f}_{L^2(m)}^2}.
  \]
The non-negative definition of $f\mapsto e^{-\lambda\norm{f}_{L^2(m)}^2}$ can be written
\begin{align*}
\sum_{i=1}^k\sum_{j=1}^k u_i u_j e^{-\lambda\norm{f_i-f_j}^2_{L^2}} 
= \sum_{i=1}^k\sum_{j=1}^k u_i u_j\espb{e^{i\sqrt{2\lambda}\scpr{X}{f_i-f_j}}} 
= \normbb{\sum_{i=1}^k u_i e^{i\sqrt{2\lambda}\scpr{X}{f_i}}}_{L^2(\Omega)}^2 \geq 0.
\end{align*}
For any $U_1,\ldots,U_k\in\Ai$, the previous result is applied to $f_1=\indi_{U_1},\ldots,f_k=\indi_{U_k}\in L^2(m)$.
As in the proof of Lemma 2.9 in \cite{Herbin.Merzbach(2006)}, we remark that
\begin{align*}
\forall i,j\in\brc{1,\dotsc,k};\quad  
m(U_i\disy U_j) 
= m\pthb{\abs{ \indi_{U_i}-\indi_{U_j} }} 
= \normb{f_i-f_j}_{L^2(m)}^2,
\end{align*}
and we deduce
\[
\sum_{i=1}^k\sum_{j=1}^k u_i u_j e^{-\lambda m(U_i\disy U_j)} \geq 0
\]
which proves the result.
\end{proof}

\noindent According to Lemma \ref{lem:cov}, we can define

\begin{definition}\label{def:ssiou}
Given the indexing collection $\mathcal{A}$ and positive real numbers $\lambda$ and $\sigma$, any Gaussian process $\{X_U;\;U\in\mathcal{A}\}$ such that for all $U,V\in\mathcal{A}$,
\begin{equation*}
\esp{X_U} = 0\quad\textrm{and}\quad \esp{X_U X_V} = \frac{\sigma^2}{2\lambda} e^{-\lambda m(U\disy V)},
\end{equation*}
is called a {\em stationary set-indexed Ornstein-Uhlenbeck (ssiOU) process}.
\end{definition}

The covariance structure of the Gaussian process coming from Definition~\ref{def:ssiou} directly implies the $L^2$-continuity and stationarity properties.

\begin{proposition} \label{prop:ssiou_continuity}
The stationary set-indexed Ornstein-Uhlenbeck process $X$ of Definition~\ref{def:ssiou} is $L^2$-monotone inner- and outer-continuous, i.e. for any increasing sequence $(U_n)_{n\in\N}$ in \Ai, such that $\overline{\cup_{k\in\N}U_k}\in\Ai$ and for any decreasing sequence $(V_n)_{n\in\N}$ in \Ai,
\[
\lim_{n\rightarrow\infty} \espb{\abs{X_{U_n}-X_{\overline{\cup_{k\in\N}U_k}}}^2} = 0
\quad\text{and}\quad 
\lim_{n\rightarrow\infty} \espb{\abs{X_{V_n}-X_{\cap_{k\in\N}V_k}}^2} = 0.
\]
\end{proposition}

\begin{proof}
  Let $(U_n)_{n\in\N}$ be an increasing sequence in \Ai such that $\overline{\cup_{k\in\N}U_k}\in\Ai$. Then using equation \eqref{eq:ssiou_cov}, we have
\begin{align*}
\forall n\in\N;\quad 
\espb{\abs{X_{U_n}-X_{\overline{\cup_{k\in\N}U_k}}}^2} 
= \frac{\sigma^2}{2\lambda} \pthb{ 2 - 2e^{-\lambda m\pth{\overline{\cup_{k\in\N}U_k}\setminus U_n}} }.
\end{align*}
  According to Assumption \textit{(iii)} on \Ai and $m$, $\lim_{n\rightarrow \infty} m\pth{\overline{\cup_{k\in\N}U_k}\setminus U_n} = 0$. Therefore, the $L^2$-monotone inner-continuity follows, and similarly, the outer-continuity of $X$.
\end{proof}

The stationarity increments property for set-indexed processes has been introduced in \cite{Herbin.Merzbach(2009)} in the context of fractional Brownian motion, and it has constitued the key property to derive deep understanding of the set-indexed L\'evy processes in \cite{Herbin.Merzbach(2011)}. The stationarity property defined below is closely related to these two previous works.

\begin{proposition} \label{prop:ssiou_stationarity}
The stationary set-indexed Ornstein-Uhlenbeck process $X$ of Definition~\ref{def:ssiou} is $m$-stationary, i.e. for any $k\in\N$, $V\in\Ai$ and increasing sequences $(U_i)_{1\leq i\leq k}$ and $(A_i)_{1\leq i\leq k}$ in \Ai such that $m(U_i\setminus V)=m(A_i)$ for all $i\in\brc{1,\dotsc,k}$, $X$ satisfies
\[
\pth{X_{U_1},\ldots,X_{U_k}} \stackrel{(d)}{=} \pth{X_{A_1},\ldots,X_{A_k}}.
\]
\end{proposition}

\begin{proof}
Let $V$, $(U_i)_{1\leq i\leq k}$ and $(A_i)_{1\leq i\leq k}$ be as in the statement. Without any loss of generality, we suppose that $V\subseteq U_i$. Then, for all $j\geq i$, as $U_i\subseteq U_j$ and $A_i\subseteq A_j$,
  \begin{align*}
    m(U_i\disy U_j)
    = m(U_j) - m(U_i)
    = m(U_j\setminus V)-m(U_i\setminus V) 
    = m(A_j) - m(A_i) 
    = m(A_j\disy A_i).
  \end{align*}
Therefore, we deduce the expected equality, since $X$ is a centered Gaussian process and for all $i,j\in\brc{1,\dotsc,k}$,
$\esp{X_{U_i}X_{U_j}}
    = \displaystyle\frac{\sigma^2}{2\lambda} e^{-\lambda m(U_i\disy U_j)} 
    = \frac{\sigma^2}{2\lambda} e^{-\lambda m(A_j\disy A_i)} 
    = \esp{X_{A_i}X_{A_j}}$.
\end{proof}

We observe that the definition of stationarity is given in a strict sense, since it concerns the invariance of finite-dimensional distributions under a form of {\it measure-invariant translation}. In the classic theory of stationary random fields, a weaker property relying on the correlation function is usually defined (see \cite{Yaglom(1987)}): $C(s,t)=\esp{X_sX_t}$ only depends on the difference $t-s$. The weak definition of stationarity for one-parameter processes can be naturally extended to the multiparameter case, but it appears that this straightforward extension is not the most relevant. Indeed, the stationarity of increments defined using Lebesgue measure or their invariance under translation appeared to be more interesting to study multiparameter processes (see e.g. L\'evy and fractional Brownian sheets), and this fact explains the form of the set-indexed extension for the stationarity property.

\subsection{Markov property and characterisation of the stationary set-indexed Ornstein-Uhlenbeck process}

To investigate the Markov property, we first need to recall a few notations used in \cite{Balanca(2012)}. Let $C\in\Ci$ such that $C = A\setminus B$, with $B\in\Aiu$ and $B\subseteq A\in\Ai$. Since the assumption \emph{Shape} holds on \Ai, Definition 1.4.5 in \cite{Ivanoff.Merzbach(2000)} states that there exists a unique \emph{extremal representation} $\brc{A_i}_{i\leq k}$ of $B$, i.e. such that $B = \cup_{i=1}^k A_i$ and for all $i\neq j$, $A_i\nsubseteq A_j$.

Then, let $\Ai_\ell$ be the semilattice $\brc{A_1\cap\dotsb\cap A_k,\dotsc,A_1\cap A_2,A_1\dotsc,A_k}\subset\Ai$ and \Cp be defined as the following subset of $\Ai_\ell$,
\begin{equation}\label{eq:def_Cfrontier}
  \Cp = \brc{U\in\Ai_\ell;\, U\nsubseteq B^\circ} \eqdef \brc{\Cpi{1},\dotsb,\Cpi{n}}, \qquad\text{where }n=\#(\Cp).
\end{equation}
The notation \vXCp{C} refers to the random vector $\vXCp{C} = \pthb{ X_{\Cpi{1}},\dotsc,X_{\Cpi{n} } }$, and similarly \vxCp{C} is used for a vector of variables. Thereby, according to \cite{Balanca(2012)}, the extension $\Delta X$ of $X$ on the class \Ci satisfies
\begin{align} \label{eq:inc_exc_formula}
\Delta X_C &\eqdef X_A - \bktbb{ \sum_{i=1}^k X_{A_i} - \sum_{i<j} X_{A_i \cap A_j} + \dotsb + (-1)^{k+1} X_{A_1\cap\dotsb\cap A_k} } \nonumber\\
&= X_A - \bktbb{ \sum_{i=1}^{n} (-1)^{\eCpi{i}} X_{\Cpi{i}} },
\end{align}
where $(-1)^{\eCpi{i}}$ represents the sign of the term $X_{\Cpi{i}}$ in the inclusion-exclusion formula. 
In other words, \eqref{eq:inc_exc_formula} says that every term $X_U$ in the previous inclusion-exclusion formula such that $U\notin\Cp$ is cancelled by another term in the sum. 

Finally $\brc{ \Fi_A; A\in\Ai }$ denotes the natural filtration generated by $X$, and for all $B\in\Aiu$ and $C\in\Ci$, $\Fi_B$ and $\Gc{C}$ respectively correspond to
\begin{equation} \label{eq:def_filtrations}
  \Fi_B = \bigvee_{A\in\Ai,A\subseteq B}\,\Fi_A \qquad\text{and}\qquad \Gc{C} = \bigvee_{B\in\Aiu,B\cap C=\vset}\,\Fi_B.
\end{equation}
We note that these filtrations are not necessarily outer-continuous.

In the following result, we prove that the ssiOU process satisfies the \Ci-Markov property introduced in \cite{Balanca(2012)}.

\begin{proposition} \label{prop:ssiou_c_markov}
The stationary set-indexed Ornstein-Uhlenbeck process $X$ of Definition~\ref{def:ssiou} is a \Ci-Markov process with respect to its natural filtration $(\Fi_A)_{A\in\Ai}$, i.e. for all $C = A\setminus B$ with $A\in\Ai$, $B\in\Ai(u)$ and all Borel function $f:\R\rightarrow\R_+$, $X$ satisfies
\begin{equation}  \label{eq:def_cmarkov}
  \espc{f(X_A)}{\Gc{C}} = \espc{f(X_A)}{\vXCp{C}} \eqdef \Pc{C} f(\vXCp{C}) \quad \Pr\textrm{-a.s}.
\end{equation}
\end{proposition}

\begin{proof}
Let $C=A\setminus B$ be in \Ci, $\brc{A_i}_{i\leq k}$ be the extremal representation of $B$ and $U$ be in $\Ai$ such that $U\cap C=\vset$. We first note that $U\cap A = (U\cap C)\cup(U\cap B) = U\cap B\in\Ai$. Thus, since \Ai satisfies the \emph{Shape} hypothesis, there exists $l\in\brc{1,\dotsc,k}$ such that $U\cap B= U\cap A_l$.
Consider the following quantity $I_U$,
\begin{align*}
I_U &= \espB{ X_U \pthB{ X_A - \sum_{i=1}^k X_{A_i}e^{-\lambda m(A\setminus A_i)} + \sum_{1\leq i<j \leq k} X_{A_i\cap A_j} e^{-\lambda m(A\setminus A_i\cap A_j)} + \\
& \qquad\qquad\qquad\qquad \dotsb + (-1)^{k} X_{A_1\cap \dotsb \cap A_k}e^{-\lambda m(A\setminus A_1\cap \dotsb \cap A_k)} } } \\
&= \frac{\sigma^2}{2\lambda} \pthB{ e^{-\lambda m(A \disy U)} - \sum_{i=1}^k e^{-\lambda (m(A_i \disy U) + m(A\setminus A_i))} + \sum_{1\leq i<j \leq k} e^{-\lambda ( m(A_i\cap A_j \disy U) + m(A\setminus A_i\cap A_j))} + \\
& \qquad\qquad\qquad\qquad \dotsb + (-1)^{k} e^{-\lambda (m(A_1\cap \dotsb \cap A_k \disy U)+m(A\setminus A_1\cap \dotsb \cap A_k))} } \\
&= \frac{\sigma^2}{2\lambda} e^{-\lambda (m(A)+m(U))} \pthB{ e^{-2\lambda m(A\cap U)} - \sum_{i=1}^k e^{-2\lambda m(A_i \cap U)} + \sum_{1\leq i<j \leq k} e^{-2\lambda m(A_i\cap A_j\cap U)} + \\
& \qquad\qquad\qquad\qquad \dotsb + (-1)^{k} e^{-2\lambda m(A_1\cap \dotsb \cap A_k \cap U)} }.
\end{align*}

\noindent Let us introduce the set-indexed function $h:A\mapsto e^{-2\lambda m(A\cap U)}$. Since the assumption \emph{Shape} holds, $h$ admits an extension $\Delta h$ on \Aiu based on an inclusion-exclusion formula. Thus, we have $I_U = \frac{\sigma^2}{2\lambda} e^{-\lambda (m(A)+m(U))} \pthb{ h(A\cap U) - \Delta h(B\cap U)}$. But since $A\cap U = B\cap U = A_l\cap U\in\Ai$ and $h$ coincides with $\Delta h$ on \Ai, we obtain $I_U = 0$. 
  
Therefore, $I_U = 0$ for all $U\in\Ai$ such that $U\cap C=\vset$ and as $X$ is a Gaussian process, we can claim that the random variable
\[
X_A - \sum_{i=1}^k X_{A_i}e^{-\lambda m(A\setminus A_i)} + \dotsb + (-1)^{k} X_{A_1\cap \dotsb \cap A_k}e^{-\lambda m(A\setminus A_1\cap \dotsb \cap A_k)} 
\]
and $\Gc{C}$ are independent.
Since the previous random variable is expressed as an inclusion-exclusion formula, equation \eqref{eq:inc_exc_formula} shows that it can also be expressed as
\[
X_A - Z_C \quad\textrm{with}\quad 
Z_C = e^{-\lambda m(A)}\bktbb{ \sum_{i=1}^{n} (-1)^{\eCpi{i}} X_{\Cpi{i}}\,e^{\lambda m(\Cpi{i})} }.
\]
Notice that $Z_C$ is $\vXCp{C}$-measurable (and then $\Gc{C}$-measurable) and $X_A - Z_C$ is independent of $\Gc{C}$ (and then independent of $\vXCp{C}$). Hence, using a classic property of the conditional expectation, we have
\[
  \espc{f(X_A)}{\Gc{C}} = \esp{f(X_A-Z_C + Z_C) \mid \Gc{C}} = \espc{f(X_A)}{Z_C}.
\]
We similarly obtain the equality $\espc{f(X_A)}{\vXCp{C}} = \espc{f(X_A)}{Z_C}$ which ends the proof.
\end{proof}
Intuitively, the \Ci-Markov property can be understood as following: For any increment $C = A\setminus B$, the $\sigma$-field $\Gc{C}$ represents the past, described as \emph{strong} as it contains all the information inside the regions $B$ satisfying $C = A\setminus B$. The vector \vXCp{C} itself gathers the minimum information related to the "border" points of $C$, and finally, $X_A$ represents the future value of the process. Then, Equation \eqref{eq:def_cmarkov} simply states that conditioning the future with respect the full history $\Gc{C}$ or the vector $\vXCp{C}$ are equivalent.

According to Proposition 2.9 in \cite{Balanca(2012)}, we can deduce that the set-indexed Ornstein-Uhlenbeck process also satisfies set-indexed \emph{sharp-Markov} and \emph{Markov} properties whose definitions can be found in \cite{Ivanoff.Merzbach(2000)a}. 
In the multiparameter case, it implies that this process is \emph{sharp-Markov} and \emph{germ-Markov} with respect to finite unions of rectangles (see \cite{Ivanoff.Merzbach(2000)a, Balanca(2012)}). The question whether this implication remains true for more complex sets has not been investigated yet (see \cite{Dalang.Walsh(1992),Dalang.Walsh(1992)a} for answers in the particular case of Brownian and L\'evy sheets).

As a consequence of the previous proposition, we can derive the \Ci-transition system \Pii and the initial law that characterize entirely a ssiOU process.

\begin{corollary} \label{cor:ssiou_markov_dens} 
The \Ci-transition system $\Pii = \brc{\Pc{C}(\vxCp{C};\Gamma);\,C\in\Ci, \Gamma\in\mathcal{B}(\R)}$ of the stationary set-indexed Ornstein Uhlenbeck process of Definition~\ref{def:ssiou} is characterized by the following transition densities, for all $C=A\setminus B\in\Ci$:
\begin{equation} \label{eq:ssiou_markov_dens}
\pc{C}(\vxCp{C};y) = \frac{1}{\sigma_C\sqrt{2\pi}} \exp\bktBB{ -\frac{1}{2\sigma_C^2} \pthbb{ y - e^{-\lambda m(A)}\bktbb{ \sum_{i=1}^{n} (-1)^{\eCpi{i}}\, x_{\Cpi{i}}\, e^{\lambda m(\Cpi{i})} } }^2 },
\end{equation}
where
\[
\sigma_C^2 = \frac{\sigma^2}{2\lambda} \pthbb{ 1 - e^{-2\lambda m(A)} \bktbb{ \sum_{i=1}^{n} (-1)^{\eCpi{i}} e^{2\lambda m(\Cpi{i})} } }.
\]
Furthermore, the initial law is given by $X_\zset\sim\Ni(0,\tfrac{\sigma^2}{2\lambda})$.
\end{corollary}

\begin{proof}
Let $C = A\setminus B$ be in \Ci and let $Z_C$ and $Y_C$ be the following Gaussian variables
\[
Z_C = e^{-\lambda m(A)}\bktbb{ \sum_{i=1}^{n} (-1)^{\eCpi{i}} X_{\Cpi{i}}\,e^{\lambda m(\Cpi{i})} } \quad\text{and}\quad Y_C = X_A - Z_C.
\]
Since the process $X$ is centered, $\esp{Z_C} = \esp{Y_C} = 0$. We note $\sigma_C^2$ the variance of $Y_C$. Using the independence of $Y_C$ and $\Gi^{*}_C$, shown in the proof of Proposition \ref{prop:ssiou_c_markov}, and the fact that $Z_C$ is $\Gi^{*}_C$-measurable, we have for any measurable function $f:\R\rightarrow\R_+$,
\begin{align*}
\espc{f(X_A)}{\Gi^{*}_C}
= &\espc{f( Z_C + Y_C )}{\Gi^{*}_C} \\
= &\frac{1}{\sigma_C\sqrt{2\pi}}  \int_\R f\pth{u + Z_C} \exp\pthbb{-\frac{ u^2}{2\sigma_C^2}} \dt u\\
= &\frac{1}{\sigma_C\sqrt{2\pi}}  \int_\R f\pth{v} \exp\pthbb{-\frac{ (v-Z_C)^2}{2\sigma_C^2}} \dt v
\eqdef \int_\R f(v)\, \pc{C}(\vXCp{C};v) \,\dt v.
\end{align*}
Equation \eqref{eq:ssiou_markov_dens} follows from this last equality. 
It remains to prove the expression of the variance $\sigma_C^2$. We first note that, as $X_{\Cpi{i}}$ is $\Gi^{*}_C$-measurable and $Y_C$ is independent of $\Gi^{*}_C$,  $\esp{X_{\Cpi{i}} Y_C} = 0$ for any $i\in\brc{1,\dotsc,n}$. Therefore, we have
\begin{align*}
\sigma_C^2 = \esp{X_A Y_C}
&= \esp{X_A^2} - e^{-\lambda m(A)} \bktbb{ \sum_{i=1}^{n} (-1)^{\eCpi{i}} \,\espb{X_A X_{\Cpi{i}}} \,e^{\lambda m(\Cpi{i})} } \\
&= \frac{\sigma^2}{2\lambda}\pthbb{ 1 - e^{-\lambda m(A)} \bktbb{ \sum_{i=1}^{n} (-1)^{\eCpi{i}} \,e^{-\lambda m(A\disy \Cpi{i})}\,e^{\lambda m(\Cpi{i})} } } \\
&= \frac{\sigma^2}{2\lambda} \pthbb{ 1 - e^{-2\lambda m(A)} \bktbb{ \sum_{i=1}^{n} (-1)^{\eCpi{i}} e^{2\lambda m(\Cpi{i})} } }.
\end{align*}
\end{proof}

The following result shows that properties exhibited in Propositions \ref{prop:ssiou_continuity}, \ref{prop:ssiou_stationarity} and  \ref{prop:ssiou_c_markov} lead to a complete characterization of the stationary set-indexed Ornstein-Uhlenbeck process.

\begin{theorem} \label{th:ssiou_charac}
A set-indexed mean-zero Gaussian process $X=\brc{X_U;\,U\in\Ai}$ is a stationary set-indexed Ornstein-Uhlenbeck process if and only if the three following properties hold:
\begin{enumerate}[(i)]
\item $L^2$-monotone inner- and outer-continuity;
\item $m$-stationarity;
\item \Ci-Markov property.
\end{enumerate}
\end{theorem}

\begin{proof}
We already know that the stationary set-indexed Ornstein-Uhlenbeck process of Definition~\ref{def:ssiou} satisfies these three properties.
Conversely, let $X$ be a zero-mean Gaussian set-indexed process which is $L^2$-monotone inner- and outer-continuous, $m$-stationary and \Ci-Markov. Without any loss of generality, we suppose $\esp{X_\zset^2}=1$.
  
We first consider an increasing and continuous function $f:\R_+\rightarrow\Ai$, i.e. an {\em elementary flow} in the terminology of \cite{Ivanoff.Merzbach(2000), Herbin.Merzbach(2009)}, such that $f(0)=\zset$. Since $m$ is monotonically continuous on \Ai (Condition {\it (iii)} of the indexing collection), the function $\theta:t\mapsto m[f(t)]$ is continuous, $\theta(0)=0$ and the pseudo-inverse $\theta^{-1}(t)=\inf\{u:\theta(u)>t\}$ satisfies $\theta\circ\theta^{-1}(t)=t$. 
Then, the projected one-parameter process $X^{m,f} = \{X_{f\circ\theta^{-1}(t)};\;t\in\R_+\}$ is a centered one-parameter Gaussian process which is $L^2$-continuous, stationary (see \cite{Herbin.Merzbach(2009)}) and Markov (see \cite{Balanca(2012)}, Proposition 2.10). 
Therefore, $X^{m,f}$ is a one-dimensional Ornstein-Uhlenbeck process (see e.g. \cite{Samorodnitsky.Taqqu(1994)}).  Since $\esp{(X^f_0)^2}=1$, there exists $\lambda_f > 0$ such that for all $s,t\in\R_+$,
\begin{align*}
\esp{X^{m,f}_s X^{m,f}_t} = e^{-\lambda_f \abs{t-s}} = e^{-\lambda_f \abs{m[f\circ\theta^{-1}(t)]-m[f\circ\theta^{-1}(s)]}} = e^{-\lambda_f m[f\circ\theta^{-1}(s)\disy f\circ\theta^{-1}(t)] }.
\end{align*}

Let us prove the constant $\lambda_f$ does not depend on the function $f$. Let $f_1$ and $f_2$ be two different elementary flows which satisfy the previous conditions. Then, as $m(\zset)=0$, for any $t>0$, we know that $m(f_1\circ\theta_1^{-1}(t)\setminus\zset) = m(f_2\circ\theta_2^{-1}(t))=t$, and therefore, according to the $m$-stationarity of $X$,
$\pthb{X^{m,f_1}_t,X^{m,f_1}_0} \eqd \pth{X^{m,f_2}_t,X^{m,f_2}_0}$,
which implies $e^{-\lambda_{f_1} t} = \esp{X^{m,f_1}_t X^{m,f_1}_0} = \esp{X^{m,f_2}_t X^{m,f_2}_0} = e^{-\lambda_{f_2} t}$, i.e. $\lambda_{f_1} = \lambda_{f_2}\eqdef\lambda$.
  
For all $U,V\in\Ai$ such that $U\subseteq V$, there exists $f$ which goes through $U$ and $V$. We obtain
\[
\esp{X_U X_V} = e^{-\lambda \abs{m(V)-m(U)}} = e^{-\lambda m(U\disy V)}.
\]
Finally let $U,V\in\Ai$. From the previous equation, we observe that 
\[
\esp{(X_V - e^{-\lambda m(V\setminus U)} X_{U\cap V} )X_{U\cap V}} = e^{-\lambda m(V\disy (U\cap V))} - e^{-\lambda m(V\setminus U)} = 0.
\]
Therefore, since $X$ is a Gaussian process, $\espc{X_V}{X_{U\cap V}} = e^{-\lambda m(V\setminus U)} X_{U\cap V}$, and using the \Ci-Markov property applied to $C=U\setminus V$ with the fact $\vXCp{U\setminus V}=X_{U\cap V}$, we obtain the expected covariance,
\begin{align*}
\esp{X_U X_V} = \espb{ X_U \,\espc{X_V}{\Gi^{*}_{U\setminus V}} } 
&= \espb{ X_U \,\espc{X_V}{X_{U\cap V}} } \\
&= e^{-\lambda m(V\setminus U)} \, \espb{ X_U  X_{U\cap V} } \\
&= e^{-\lambda m(V\setminus U)} \cdot e^{-\lambda m(U\setminus V)} = e^{-\lambda m(U\disy V)}.
\end{align*}
\end{proof}

\section{Definition of a general \emph{set-indexed Ornstein-Uhlenbeck process}}
\label{sec:gensiou}

Using the \Ci-Markov property obtained in Proposition \ref{prop:ssiou_c_markov} and the \Ci-transition system \Pii from Corollary \ref{cor:ssiou_markov_dens}, we can finally define a general \emph{set-indexed Ornstein-Uhlenbeck process}.

\begin{definition} \label{def:siou}
A process $X$ is called a \emph{set-indexed Ornstein-Uhlenbeck process} if
\begin{enumerate}[\it (i)]
\item $X_\zset\sim\nu$, where $\nu$ is a given initial probability distribution;
\item $X$ is \Ci-Markov with a \Ci-transition system given by \eqref{eq:ssiou_markov_dens}.
\end{enumerate}
\end{definition}

Theorem 2.2 in \cite{Balanca(2012)} proves the existence of such processes in the canonical probability space $(\R^{\Ai}, \Pr_{\nu})$ for any initial probability distribution $\nu$. Then, $\Pr_{\nu}$ is the probability measure on $\R^{\Ai}$ under which the canonical process defined by $X_U(\omega)=\omega(U)$ for all $\omega\in\R^{\Ai}$ is a set-indexed Ornstein-Uhlenbeck process.
In the particular case of Dirac initial distribution, the complete determination of the laws of $X$ is given by the following result.

\begin{proposition} \label{prop:siou_deltax}
For any $x\in\R$, under the probability $\Pr_x$, the canonical set-indexed Ornstein-Uhlenbeck process $X$ is the Gaussian process defined by the covariance structure
\begin{align} \label{eq:siou_deltax1}
\forall U\in\Ai;\quad &\esp[_x]{X_U} = x\ e^{-\lambda m(U)}, \\
\forall U,V\in\Ai;\quad &\cov_x(X_U,X_V) = \frac{\sigma^2}{2\lambda}\pthb{ e^{-\lambda m(U\disy V)} - e^{-\lambda( m(U)+m(V) )} }.\label{eq:siou_deltax2}
\end{align}
\end{proposition}

\begin{proof}
We first check that $X$ is a Gaussian process under the probability $\Pr_x$. \\
Let $A_1,\dotsc,A_k\in\Ai$ and $\lambda_1,\dotsc,\lambda_k\in\R$. Without any loss of generality, we can suppose that $\Ai_\ell = \brc{A_0=\zset,A_1,\dotsc,A_k}$ is a semilattice and we denote $C_i = A_i\setminus(\cup_{j=0}^{i-1} A_j)$ for all $i\in\brc{1,\dotsc,k}$. Then, using notations from Corollary \ref{cor:ssiou_markov_dens}, we have
\begin{align*}
\espbb[_x]{\exp\pthbb{ i\sum_{j=1}^k \lambda_j X_{A_j} }}
&= \espbb[_x]{\exp\pthbb{ i\sum_{j=1}^{k-1} \lambda_j X_{A_j} } \,\espcB[_x]{\exp\pthb{i\lambda_k X_{A_k}}}{\Gi_{C_k}^*} } \\
&= \espbb[_x]{\exp\pthbb{ i\sum_{j=1}^{k-1} \lambda_j X_{A_j} } \,\exp\pthb{i\lambda_k Z_{C_k}} \,\espb[_x]{\exp\pth{i\lambda_k Y_{C_k}}} } \\
&= \exp\pthbb{-\frac{\lambda_k^2\sigma^2_{C_k}}{2}} \,\espbb[_x]{\exp\pthbb{ i\sum_{j=1}^{k-1} \lambda'_j X_{A_j} } },
\end{align*}
since $Z_{C_k}$ is weighted sum of $X_V$, $V\in\brc{A_0,\dotsc,A_{k-1}}$.
Therefore, by induction on $k$, we get the characteristic function of a Gaussian variable.
  
In order to obtain the mean and the covariance functions, we consider the case $k=3$, with the semi-lattice $\{\zset, A_1=A_2\cap A_3, A_2, A_3\}$.
We compute
\begin{align*}
\espB[_x]{\exp\pthB{ i( \lambda_2 X_{A_2} + \lambda_3 X_{A_3} ) } } 
&= \exp\pthB{-\frac{1}{2}\lambda_3^2\sigma^2_{C_3}} \,\espB[_x]{\exp\pthB{ i( \lambda_2 X_{A_2} + \lambda_3 Z_{C_3} ) } } \\
&= \exp\pthB{-\frac{1}{2}\lambda_3^2\sigma^2_{C_3}} \,\espB[_x]{\exp\pthB{ i( \lambda_2 X_{A_2} + \lambda_3 e^{-\lambda m(A_3\setminus A_1)} X_{A_1} ) } }.
\end{align*}
Using the \Ci-Markov property applied to $C_2=A_2\setminus A_1$, we get
\begin{align*}
\espB[_x]{\exp\pthB{ i( \lambda_2 X_{A_2} + \lambda_3 X_{A_3} ) } } 
&= \exp\pthB{-\frac{1}{2}\lambda_3^2\sigma^2_{C_3} - \frac{1}{2}\lambda_3^2\sigma^2_{C_2}} \\
&\qquad\times\espB[_x]{\exp\pthB{ i\pthb{ \lambda_2 e^{-\lambda m(A_2\setminus A_1)} + \lambda_3 e^{-\lambda m(A_3\setminus A_1)} } X_{A_1} } }.
\end{align*}
Then, the \Ci-Markov property applied to $C_1=A_1\setminus\zset$ leads to
\begin{align*}
&\espB[_x]{\exp\pthB{ i( \lambda_2 X_{A_2} + \lambda_3 X_{A_3} ) } } \\
&\qquad\qquad = \exp\pthB{-\frac{1}{2}\lambda_3^2\sigma^2_{C_3} - \frac{1}{2}\lambda_3^2\sigma^2_{C_2} - \frac{1}{2}\pthb{ \lambda_2 e^{-\lambda m(A_2\setminus A_1)} + \lambda_3 e^{-\lambda m(A_3\setminus A_1)} }^2 \sigma_{C_1}^2 } \\
&\qquad\qquad\qquad\qquad\times\espB[_x]{\exp\pthB{ i \pthb{ \lambda_2 e^{-\lambda m(A_2\setminus A_1)} + \lambda_3 e^{-\lambda m(A_3\setminus A_1)}}  e^{-\lambda m(A_1)} X_\zset } } \\
&\qquad\qquad = \exp\pthB{-\frac{1}{2}\lambda_3^2\sigma^2_{C_3} - \frac{1}{2}\lambda_3^2\sigma^2_{C_2} - \frac{1}{2}\pthb{ \lambda_2 e^{-\lambda m(A_2\setminus A_1)} + \lambda_3 e^{-\lambda m(A_3\setminus A_1)} }^2 \sigma_{C_1}^2 } \\
&\qquad\qquad\qquad\qquad\times\exp\pthB{ i\pthb{ \lambda_2 e^{-\lambda m(A_2)} + \lambda_3 e^{-\lambda m(A_3)} } x }.
\end{align*}
The mean of $X$ comes from the last line. The covariance is obtained from the cross term in front of $\lambda_2\lambda_3$:
  \begin{align*}
    \sigma_{C_1}^2 e^{-\lambda m(A_2\setminus A_1)} e^{-\lambda m(A_3\setminus A_1)} 
    &= \frac{\sigma^2}{2\lambda}\pthb{1 - e^{-2\lambda m(A_1)}} e^{-\lambda m(A_2\disy A_3)} \\
    &= \frac{\sigma^2}{2\lambda}\pthb{ e^{-\lambda m(A_2\disy A_3)} - e^{-\lambda( m(A_2)+m(A_3) )} },
  \end{align*}
  since $A_1=A_2\cap A_3$ and $\sigma_{C_1}^2 = \frac{\sigma^2}{2\lambda}\pthb{1 - e^{-2\lambda m(A_1)}}$.
\end{proof}

\section{Multiparameter Ornstein-Uhlenbeck process} \label{sec:mpou}

In the particular case of the indexing collection $\Ai=\brc{\ivff{0,t}; t\in\R_+}$ endowed with the Lebesgue measure $m$, the set-indexed Ornstein-Uhlenbeck processes studied in Sections \ref{sec:siou} and \ref{sec:gensiou} reduce to the classical one-dimensional Ornstein-Uhlenbeck process. 

In the multiparameter setting, a natural extension of the stationary Ornstein-Uhlenbeck process can be defined by
\begin{equation} \label{eq:mpsou}
  \forall t\in\R_+^N;\quad Y_t = \int_{-\infty}^{t} \sigma\, e^{\scpr{\alpha}{u-t}}\dt W_u.
\end{equation}
where $\sigma > 0$, $\alpha=(\alpha_1,\dots,\alpha_N)\in\R^N$ with $\alpha_i > 0$ and $W$ is the Brownian sheet. The covariance of this process is given by
\[
  \esp{Y_s Y_t} = \prod_{i=1}^N \int_{-\infty}^{s_i\wedge t_i} \sigma^2 e^{\alpha_i(u_i-s_i-t_i)}\,\dt u_i = \frac{\sigma^2}{\prod_{i=1}^N \alpha_i} \exp\brcbb{ -\sum_{i=1}^N \alpha_i (s_i+t_i - s_i\wedge t_i) }.
\]
Hence, $Y$ is a stationary set-indexed Ornstein-Uhlenbeck process on the space $\R_+^N$ endowed with the indexing collection $\Ai=\brc{\ivff{0,t}; t\in\R^N_+}$ and the measure $m_\alpha$ defined on the Borel $\sigma$-field by
\begin{equation}\label{eq:m-alpha}
\forall A\in\Bi(\R^N); \quad m_\alpha(A) = \sum_{i=1}^N \alpha_i \lambda_1 \pth{A\cap e_i},
\end{equation}
where $\lambda_1$ is the Lebesgue measure on $\R$ and $e_1,\dotsc,e_N$ are the axes of $\R^N$: $e_1=\R\times\brc{0}^{N-1}$, $e_2=\brc{0}\times\R\times\brc{0}^{N-2}$, \ldots

\medskip

The following proposition extends this result to the general set-indexed Ornstein-Uhlenbeck process defined in Section \ref{sec:gensiou}, proving that it also has a natural integral representation in the particular multiparameter case.

\begin{proposition} \label{prop:siou_mpou}
  Let $Y=\{Y_t;\;t\in\R_+^N\}$ be the multiparameter process defined by
  \begin{equation} \label{eq:mpou}
    \forall t\in\R_+^N;\quad Y_t = e^{-\scpr{\alpha}{t}}\bktbb{Y_0 + \sigma\int_{\ivof{-\infty,t}\setminus\ivof{-\infty,0}} e^{\scpr{\alpha}{u}}\dt W_u },
  \end{equation}
  where $\sigma > 0$, $\alpha=(\alpha_1,\dots,\alpha_N)\in\R^N$ with $\alpha_i > 0$ for each $i\in\brc{1,\dotsc,N}$, $W$ is the Brownian sheet and $Y_0$ is a random variable independent of $W$.

  \noindent Then, $Y$ is a set-indexed Ornstein-Uhlenbeck process of Definition~\ref{def:siou} on the space $(\Ti,\Ai,m_\alpha)$, with $\Ai=\{[0,t];\;t\in\R_+^N\}$ and $m_{\alpha}$ defined in \eqref{eq:m-alpha}.
\end{proposition}
\begin{proof}
  First we observe that the measure $m_\alpha$ satisfies, for all $s,t\in\R_+^N$,
  \[
    m_\alpha\pthb{ \ivff{0,s}\cap\ivff{0,t} } = \sum_{i=1}^N \alpha_i (s_i\wedge t_i) = \scpr{\alpha}{s\curlywedge t}\quad \text{where} \quad s\curlywedge t := (s_1\wedge t_1,\dotsc,s_N\wedge t_N).
  \]
  Let $t_1,\dotsc,t_k$ be in $\R_+^N$ and $\lambda_1,\dots,\lambda_k$ in $\R$. 
  For any fixed $x_0\in\R$, $Y^{x_0}$ denotes the Gaussian process defined by
  \[
    \forall t\in\R_+^N;\quad Y^{x_0}_t = e^{-\scpr{\alpha}{t}}\bktbb{x_0 + \sigma\int_{A_t}e^{\scpr{\alpha}{u}}\dt W_u },
  \]
  where $A_t = \ivof{-\infty,t}\setminus\ivof{-\infty,0}$.
    
  Let $Y$ be the $\R_+^N$-indexed process defined by \eqref{eq:mpou} and denote by $\nu$ the law of $Y_0$. Since $Y_0$ and $W$ are independent, we have
  \begin{align*}
    \espB{e^{i \sum_{j=1}^k \lambda_j Y_{t_j} }} 
    = \int_\R \espB{e^{i \sum_{j=1}^k \lambda_j Y^{x_0}_{t_j} }} \, \nu(\dt x_0),
  \end{align*}
  Let us determine the mean and covariance of the process $Y^{x_0}$, for any $x_0\in\R$:
  \[
  \forall t\in\R_+^N; \quad \espb{Y^{x_0}_t} = x_0 \, e^{-\scpr{\alpha}{t}} = x_0\, e^{-m_\alpha(\ivff{0,t})}
  \]
  and for all $s,t\in\R_+^N$,
  \begin{align*}
    \cov(Y^{x_0}_s,Y^{x_0}_t)
    &= \sigma^2 e^{ -\scpr{\alpha}{s+t} }\espbb{ \int_{A_s} e^{\scpr{\alpha}{u}} \dt W_u \int_{A_t} e^{\scpr{\alpha}{u}} \dt W_u }
    = \sigma^2 e^{ -\scpr{\alpha}{s+t} } \int_{A_{s\curlywedge t}} e^{2\scpr{\alpha}{u}} \dt u \\
    &= \frac{\sigma^2}{2^n \prod_{j=1}^n \alpha_j }  e^{ -\scpr{\alpha}{s+t} }\pthb{ e^{2\scpr{\alpha}{s\curlywedge t}} - 1 } \\
    &= \frac{\sigma^2}{2^n \prod_{j=1}^n \alpha_j }  e^{ -m_\alpha(\ivff{0,s})-m_\alpha(\ivff{0,t}) }\pthb{ e^{2 m_\alpha(\ivff{0,s}\cap\ivff{0,t}) } - 1 } \\
    &= \frac{\widetilde{\sigma}^2}{2} \pthB{ e^{-m_\alpha(\ivff{0,s}\disy\ivff{0,t})} -  e^{ -\pth{m_\alpha(\ivff{0,s})+m_\alpha(\ivff{0,t})} } } 
    = \cov_{x_0}(X_{\ivff{0,s}},X_{\ivff{0,t}}),
  \end{align*}
  where $X$ is the canonical set-indexed Ornstein-Uhlenbeck process with parameters $(\widetilde{\sigma},\lambda=1)$ with notations of Proposition~\ref{prop:siou_deltax}. Therefore, the process $Y^{x_0}$ has the same law as $X_{[0,\bullet]}$ starting from $x_0$ and
  \[
    \espB{e^{i \sum_{j=1}^k \lambda_j Y^{x_0}_{t_j} }} = \espB[_{x_0}]{e^{i \sum_{j=1}^k \lambda_j X_{\ivff{0,t_j}} }}.
  \]
  Consequently
  \[
    \espB{e^{i \sum_{j=1}^k \lambda_i Y_{t_j} }} = \int_\R \espB[_{x_0}]{e^{i \sum_{j=1}^k \lambda_j X_{\ivff{0,t_j}} }} \, \nu(\dt x_0) = \espB[_\nu]{e^{i \sum_{j=1}^k \lambda_j X_{\ivff{0,t_j}} }},
  \]
  which states that $Y$ and $X_{[0,\bullet]}$ have the same law and concludes the proof.
\end{proof}

\begin{remark}
We have exhibited an unusual measure $m_\alpha$ on $\R^N$, which only charges the axes $(e_i)_{i\leq N}$. This measure is also interesting when the set-indexed Brownian motion (siBM) is considered on the space $(\Ti,\Ai,m_\alpha)$ with $\alpha=(1,\dotsc,1)$, as it corresponds to a classic multiparameter process called the additive Brownian motion (see e.g. \cite{Khoshnevisan(2002)}). Conversely, since we know that the Brownian sheet is a siBM on the space $(\Ti,\Ai,\lambda)$, where $\lambda$ is the Lebesgue measure, we could also define a different multiparameter Ornstein-Uhlenbeck process using the Lebesgue measure instead of $m_\alpha$.
\end{remark}

\begin{remark}
  A different multiparameter extension of the Ornstein-Uhlenbeck process has already been introduced in the literature (e.g. see \cite{Wang(1988), Wang(1995)} and \cite{Graversen.Pedersen(2011)}). It admits an integral representation given by,
  \begin{equation} \label{eq:mpou2}
    \forall t\in\R_+^N;\quad Y_t = e^{-\scpr{\alpha}{t}}\bktbb{Y_0 + \sigma\int_0^{t}e^{\scpr{\alpha}{u}}\dt W_u }.
  \end{equation}
  If we consider a Markov point of view, the definition given in Proposition \ref{prop:siou_mpou} seems more natural. Indeed, as described in \cite{Balanca(2012)}, the transition probabilities of the process described in Equation \eqref{eq:mpou2} do not strictly correspond to those of the set-indexed Ornstein-Uhlenbeck, and can not be extended to the set-indexed formalism. Furthermore, we observe that the model \eqref{eq:mpou2} does not embrace the natural stationary case described in equation \eqref{eq:mpsou}. 
\end{remark}



\

\end{document}